\catcode`\@=11

\font\tenmsa=msam10
\font\sevenmsa=msam7
\font\fivemsa=msam5
\font\tenmsb=msbm10
\font\sevenmsb=msbm7
\font\fivemsb=msbm5
\newfam\msafam
\newfam\msbfam
\textfont\msafam=\tenmsa  \scriptfont\msafam=\sevenmsa
  \scriptscriptfont\msafam=\fivemsa
\textfont\msbfam=\tenmsb  \scriptfont\msbfam=\sevenmsb
  \scriptscriptfont\msbfam=\fivemsb

\def\hexnumber@#1{\ifnum#1<10 \number#1\else
 \ifnum#1=10 A\else\ifnum#1=11 B\else\ifnum#1=12 C\else
 \ifnum#1=13 D\else\ifnum#1=14 E\else\ifnum#1=15 F\fi\fi\fi\fi\fi\fi\fi}

\def\msa@{\hexnumber@\msafam}
\def\msb@{\hexnumber@\msbfam}
\mathchardef\boxdot="2\msa@00
\mathchardef\boxplus="2\msa@01
\mathchardef\boxtimes="2\msa@02
\mathchardef\square="0\msa@03
\mathchardef\blacksquare="0\msa@04
\mathchardef\centerdot="2\msa@05
\mathchardef\lozenge="0\msa@06
\mathchardef\blacklozenge="0\msa@07
\mathchardef\circlearrowright="3\msa@08
\mathchardef\circlearrowleft="3\msa@09
\mathchardef\rightleftharpoons="3\msa@0A
\mathchardef\leftrightharpoons="3\msa@0B
\mathchardef\boxminus="2\msa@0C
\mathchardef\Vdash="3\msa@0D
\mathchardef\Vvdash="3\msa@0E
\mathchardef\vDash="3\msa@0F
\mathchardef\twoheadrightarrow="3\msa@10
\mathchardef\twoheadleftarrow="3\msa@11
\mathchardef\leftleftarrows="3\msa@12
\mathchardef\rightrightarrows="3\msa@13
\mathchardef\upuparrows="3\msa@14
\mathchardef\downdownarrows="3\msa@15
\mathchardef\upharpoonright="3\msa@16

\mathchardef\downharpoonright="3\msa@17
\mathchardef\upharpoonleft="3\msa@18
\mathchardef\downharpoonleft="3\msa@19
\mathchardef\rightarrowtail="3\msa@1A
\mathchardef\leftarrowtail="3\msa@1B
\mathchardef\leftrightarrows="3\msa@1C
\mathchardef\rightleftarrows="3\msa@1D
\mathchardef\Lsh="3\msa@1E
\mathchardef\Rsh="3\msa@1F
\mathchardef\rightsquigarrow="3\msa@20
\mathchardef\leftrightsquigarrow="3\msa@21
\mathchardef\looparrowleft="3\msa@22
\mathchardef\looparrowright="3\msa@23
\mathchardef\circeq="3\msa@24
\mathchardef\succsim="3\msa@25
\mathchardef\gtrsim="3\msa@26
\mathchardef\gtrapprox="3\msa@27
\mathchardef\multimap="3\msa@28
\mathchardef\therefore="3\msa@29
\mathchardef\because="3\msa@2A
\mathchardef\doteqdot="3\msa@2B

\mathchardef\triangleq="3\msa@2C
\mathchardef\precsim="3\msa@2D
\mathchardef\lesssim="3\msa@2E
\mathchardef\lessapprox="3\msa@2F
\mathchardef\eqslantless="3\msa@30
\mathchardef\eqslantgtr="3\msa@31
\mathchardef\curlyeqprec="3\msa@32
\mathchardef\curlyeqsucc="3\msa@33
\mathchardef\preccurlyeq="3\msa@34
\mathchardef\leqq="3\msa@35
\mathchardef\leqslant="3\msa@36
\mathchardef\lessgtr="3\msa@37
\mathchardef\backprime="0\msa@38
\mathchardef\risingdotseq="3\msa@3A
\mathchardef\fallingdotseq="3\msa@3B
\mathchardef\succcurlyeq="3\msa@3C
\mathchardef\geqq="3\msa@3D
\mathchardef\geqslant="3\msa@3E
\mathchardef\gtrless="3\msa@3F
\mathchardef\sqsubset="3\msa@40
\mathchardef\sqsupset="3\msa@41
\mathchardef\trianglerighteq="3\msa@44
\mathchardef\trianglelefteq="3\msa@45
\mathchardef\bigstar="0\msa@46
\mathchardef\between="3\msa@47
\mathchardef\blacktriangledown="0\msa@48
\mathchardef\blacktriangleright="3\msa@49
\mathchardef\blacktriangleleft="3\msa@4A
\mathchardef\blacktriangle="0\msa@4E
\mathchardef\triangledown="0\msa@4F
\mathchardef\eqcirc="3\msa@50
\mathchardef\lesseqgtr="3\msa@51
\mathchardef\gtreqless="3\msa@52
\mathchardef\lesseqqgtr="3\msa@53
\mathchardef\gtreqqless="3\msa@54
\mathchardef\Rrightarrow="3\msa@56
\mathchardef\Lleftarrow="3\msa@57
\mathchardef\veebar="2\msa@59
\mathchardef\barwedge="2\msa@5A
\mathchardef\doublebarwedge="2\msa@5B
\mathchardef\angle="0\msa@5C
\mathchardef\measuredangle="0\msa@5D
\mathchardef\sphericalangle="0\msa@5E
\mathchardef\varpropto="3\msa@5F
\mathchardef\smallsmile="3\msa@60
\mathchardef\smallfrown="3\msa@61
\mathchardef\Subset="3\msa@62
\mathchardef\Supset="3\msa@63
\mathchardef\Cup="2\msa@64

\mathchardef\Cap="2\msa@65

\mathchardef\curlywedge="2\msa@66
\mathchardef\curlyvee="2\msa@67
\mathchardef\leftthreetimes="2\msa@68
\mathchardef\rightthreetimes="2\msa@69
\mathchardef\subseteqq="3\msa@6A
\mathchardef\supseteqq="3\msa@6B
\mathchardef\bumpeq="3\msa@6C
\mathchardef\Bumpeq="3\msa@6D
\mathchardef\lll="3\msa@6E

\mathchardef\ggg="3\msa@6F

\mathchardef\circledS="0\msa@73
\mathchardef\pitchfork="3\msa@74
\mathchardef\dotplus="2\msa@75
\mathchardef\backsim="3\msa@76
\mathchardef\backsimeq="3\msa@77
\mathchardef\complement="0\msa@7B
\mathchardef\intercal="2\msa@7C
\mathchardef\circledcirc="2\msa@7D
\mathchardef\circledast="2\msa@7E
\mathchardef\circleddash="2\msa@7F
\def\ulcorner{\delimiter"4\msa@70\msa@70 }
\def\urcorner{\delimiter"5\msa@71\msa@71 }
\def\llcorner{\delimiter"4\msa@78\msa@78 }
\def\lrcorner{\delimiter"5\msa@79\msa@79 }
\def\yen{\mathhexbox\msa@55 }
\def\checkmark{\mathhexbox\msa@58 }
\def\circledR{\mathhexbox\msa@72 }
\def\maltese{\mathhexbox\msa@7A }
\mathchardef\lvertneqq="3\msb@00
\mathchardef\gvertneqq="3\msb@01
\mathchardef\nleq="3\msb@02
\mathchardef\ngeq="3\msb@03
\mathchardef\nless="3\msb@04
\mathchardef\ngtr="3\msb@05
\mathchardef\nprec="3\msb@06
\mathchardef\nsucc="3\msb@07
\mathchardef\lneqq="3\msb@08
\mathchardef\gneqq="3\msb@09
\mathchardef\nleqslant="3\msb@0A
\mathchardef\ngeqslant="3\msb@0B
\mathchardef\lneq="3\msb@0C
\mathchardef\gneq="3\msb@0D
\mathchardef\npreceq="3\msb@0E
\mathchardef\nsucceq="3\msb@0F
\mathchardef\precnsim="3\msb@10
\mathchardef\succnsim="3\msb@11
\mathchardef\lnsim="3\msb@12
\mathchardef\gnsim="3\msb@13
\mathchardef\nleqq="3\msb@14
\mathchardef\ngeqq="3\msb@15
\mathchardef\precneqq="3\msb@16
\mathchardef\succneqq="3\msb@17
\mathchardef\precnapprox="3\msb@18
\mathchardef\succnapprox="3\msb@19
\mathchardef\lnapprox="3\msb@1A
\mathchardef\gnapprox="3\msb@1B
\mathchardef\nsim="3\msb@1C
\mathchardef\napprox="3\msb@1D
\mathchardef\nsubseteqq="3\msb@22
\mathchardef\nsupseteqq="3\msb@23
\mathchardef\subsetneqq="3\msb@24
\mathchardef\supsetneqq="3\msb@25
\mathchardef\subsetneq="3\msb@28
\mathchardef\supsetneq="3\msb@29
\mathchardef\nsubseteq="3\msb@2A
\mathchardef\nsupseteq="3\msb@2B
\mathchardef\nparallel="3\msb@2C
\mathchardef\nmid="3\msb@2D
\mathchardef\nshortmid="3\msb@2E
\mathchardef\nshortparallel="3\msb@2F
\mathchardef\nvdash="3\msb@30
\mathchardef\nVdash="3\msb@31
\mathchardef\nvDash="3\msb@32
\mathchardef\nVDash="3\msb@33
\mathchardef\ntrianglerighteq="3\msb@34
\mathchardef\ntrianglelefteq="3\msb@35
\mathchardef\ntriangleleft="3\msb@36
\mathchardef\ntriangleright="3\msb@37
\mathchardef\nleftarrow="3\msb@38
\mathchardef\nrightarrow="3\msb@39
\mathchardef\nLeftarrow="3\msb@3A
\mathchardef\nRightarrow="3\msb@3B
\mathchardef\nLeftrightarrow="3\msb@3C
\mathchardef\nleftrightarrow="3\msb@3D
\mathchardef\divideontimes="2\msb@3E
\mathchardef\varnothing="0\msb@3F
\mathchardef\nexists="0\msb@40
\mathchardef\mho="0\msb@66
\mathchardef\thorn="0\msb@67
\mathchardef\beth="0\msb@69
\mathchardef\gimel="0\msb@6A
\mathchardef\daleth="0\msb@6B
\mathchardef\lessdot="3\msb@6C
\mathchardef\gtrdot="3\msb@6D
\mathchardef\ltimes="2\msb@6E
\mathchardef\rtimes="2\msb@6F
\mathchardef\shortmid="3\msb@70
\mathchardef\shortparallel="3\msb@71
\mathchardef\smallsetminus="2\msb@72
\mathchardef\thicksim="3\msb@73
\mathchardef\thickapprox="3\msb@74
\mathchardef\approxeq="3\msb@75
\mathchardef\succapprox="3\msb@76
\mathchardef\precapprox="3\msb@77
\mathchardef\curvearrowleft="3\msb@78
\mathchardef\curvearrowright="3\msb@79
\mathchardef\digamma="0\msb@7A
\mathchardef\varkappa="0\msb@7B
\mathchardef\hslash="0\msb@7D
\mathchardef\hbar="0\msb@7E
\mathchardef\backepsilon="3\msb@7F
\def\Bbb{\ifmmode\let\next\Bbb@\else
 \def\next{\errmessage{Use \string\Bbb\space only in math mode}}\fi\next}
\def\Bbb@#1{{\Bbb@@{#1}}}
\def\Bbb@@#1{\fam\msbfam#1}

\catcode`\@=\active

\def\inv{^{\raise.15ex\hbox{${
  \scriptscriptstyle -}$}\kern-.05em 1}}

\def\Dsl{\,\raise.15ex\hbox{$/$}\mkern-13.5mu D}
\def\dsl{\raise.15ex\hbox{$/$}\kern-.57em\hbox{$\partial$}}

\def\lspace{\ifx\answ\bigans{}\else\qquad\fi}
\def\del{\partial}
 \def\CH{\hbox{{$\cal H$}}}
 
 \def\CS{\hbox{{$\cal S$}}}

\def\CR{\hbox{{$\cal R$}}}

\def\cu{\hbox{{\sl u}}} % used for Lie algebra 'm'
\def\cv{\hbox{{\sl v}}}
\def\lform{\hbox{$\sqcup$}\llap{\hbox{$\sqcap$}}}
\def\darr#1{\raise1.5ex\hbox{$\leftrightarrow$}
\mkern-16.5mu #1}
\def\h{{{1\over2}}}

\def\INT{{\textstyle \int\kern-.642em\int}}

\def\C{{\Bbb C}}

\def\eps{{\epsilon}}

\def\lcross{{>\!\!\!\triangleleft}}

\def\codcross{{\blacktriangleright\!\!\blacktriangleleft}}
\def\rbiprod{{\cdot\kern-.33em\triangleright\!\!\!<}}
\def\lbiprod{{>\!\!\!\triangleleft\kern-.33em\cdot}}

\def\tens{\mathop{\otimes}}

\def\la{{\triangleright}}

\def\ev{{\rm ev}}

\def\id{{\rm id}}

\def\<{\langle}
\def\>{\rangle}
\def\dila{{\varsigma}}

\def\haj#1{{\mathaccent20 {#1}}}
\def\Vhaj{{V\haj{\ }}}

\def\vect{{\bf t}}
\def\vecx{{\bf x}}\def\vecp{{\bf p}}
\def\vecl{{\bf l}}

\def\vecm{{\bf m}}\def\vecc{{\bf c}}

\def\<{\langle}
\def\>{\rangle}

\def\equad{\kern -1.7em}

\def\o{{}_{\scriptscriptstyle(1)}}
\def\t{{}_{\scriptscriptstyle(2)}}

\def\und#1{{\underline {#1}}}

\def\uo{{{}^{\scriptscriptstyle(1)}}}
\def\ut{{{}^{\scriptscriptstyle(2)}}}

\def\umo{{{}^{\scriptscriptstyle-(1)}}}
\def\umt{{{}^{\scriptscriptstyle-(2)}}}

\def\Bo{{{}_{\und{\scriptscriptstyle(1)}}}}
\def\Bt{{{}_{\und{\scriptscriptstyle(2)}}}}

\def\opBo{{{}_{\und{\scriptscriptstyle(1)}^{\rm op}}}}
\def\opBt{{{}_{\und{\scriptscriptstyle(2)}^{\rm op}}}}

\def\text#1{\mbox{\rm #1}}
\def\note#1{}

\def\blacksquare{{\lform}}%AMS Tex Fakes
\def\frac#1#2{{{#1\over#2}}}

\def\proof{\goodbreak\noindent{\bf Proof\quad}}

\def\endproof{{\ $\lform$}\bigskip }

\def\eqn#1#2{\begin{equation}#2\label{#1}\end{equation}}
\def\align#1{\begin{eqnarray*}#1\end{eqnarray*}}
\def\cmath#1{\[\begin{array}{c} #1 \end{array}\]}
\def\ceqn#1#2{\begin{equation}\label{#1}\begin{array}{c}#2\end{array}
\end{equation}}

\documentstyle[11pt]{article}
\newtheorem{lemma}{Lemma}[section] \newtheorem{propos}[lemma]{Proposition}
\newtheorem{example}[lemma]{Example} \newtheorem{theorem}[lemma]{Theorem}

\begin{document}\baselineskip 22pt

{\ }\qquad\qquad \hskip 4.3in  DAMTP/95-68
\vspace{.2in}

\begin{center} {\LARGE BRAIDED GEOMETRY OF  THE CONFORMAL ALGEBRA}
\\ \baselineskip 13pt{\ }
{\ }\\
Shahn Majid\footnote{Royal Society University Research Fellow and Fellow of
Pembroke College, Cambridge}\\
{\ }\\
Department of Applied Mathematics \& Theoretical Physics\\
University of Cambridge, Cambridge CB3 9EW\\
+\\
Department of Mathematics, Harvard University\\
Science Center, Cambridge MA 02138, USA\footnote{During 1995+1996}\\
\end{center}
\begin{center}
December 1995 -- revised January 1996
\end{center}

\vspace{10pt}
\begin{quote}\baselineskip 13pt
\noindent{\bf Abstract}
We show that the action of the special conformal
transformations of the usual (undeformed) conformal group is the $q\to 1$
scaling limit of the braided adjoint action or $R$-commutator of $q$-Minkowski
space on itself. We also describe the $q$-deformed conformal algebra in
$R$-matrix form and its quasi-$*$ structure.


\end{quote}
\baselineskip 22pt

\section{Introduction} It is a standard geometrical fact that the action of the
momentum generators of the Poincar\'e algebra in physics is determined by the
additive
group structure of spacetime, by differentiation as an infinitesimal addition.
In this note we provide a novel geometrical picture of the special conformal
transformations $c_i$ similarly in terms of the structure of spacetime itself.
Namely, we show that they act as the $q\to 1$ limit of the braided adjoint
action by which any braided group acts upon itself. In the case of $q$-deformed
spacetime\cite{Ma:varen}, this is the $R$-commutator
\[ c_i.(x_{i_1}\cdots x_{i_n})={x_{i_1}\cdots x_{i_n}x_i-x_{a_1}x_{j_1}\cdots
x_{j_n}{\bf R}^{j_1}{}_{i_1}{}^{a_1}{}_{a_2} {\bf
R}^{j_2}{}_{i_2}{}^{a_2}{}_{a_3}\cdots {\bf R}^{j_n}{}_{i_n}{}^{a_n}{}_i\over
q-q^{-1}}\]
where the $x_i$ are the non-commuting spacetime-coordinates with braid
statistics controlled by the appropriate Yang-Baxter matrix ${\bf
R}^i{}_j{}^k{}_l$.
We take the limit $q\to 1$. The formula also works for $q\ne 1$ and extends the
$q$-Poincar\'e action in \cite{Ma:poi} to an action of the $q$-deformed
conformal algebra.

We believe that this result is interesting as an application of $q$-deformation
techniques to classical geometry: the picture which it provides is rather
simpler than the usual picture of the $c_i$ in terms of conjugation under
conformal inversion of spacetime translation, but is only possible when $q\ne
1$. One may work with the $c_i$ in our $q\ne 1$ setting and afterwards set
$q=1$. The result adds weight to the idea that $q$-deformed geometry is
conceptually simpler and more regular than classical geometry, with certain
notions unified in ways that are singular when $q=1$.

In a previous paper\cite{Ma:lie} we  showed that the braided adjoint action
with respect to the {\em multiplicative} braided group structure\cite{Ma:exa}
of $q$-Minkowski space (as hermitian $2\times 2$ matrices being multiplied)
corresponds in  a $q\to 1$ scaling limit to the internal symmetry Lie algebra
$su(2)\times u(1)$. Our result now is for the {\em additive} braided group
structure due to Meyer\cite{Mey:new}. It appears that both adjoint actions have
important scaling limits as $q\to 1$.

We  note that while a lot of effort has been expended in $q$-deforming
geometrical structures, in particular in the author's `braided groups
approach'\cite{Ma:exa}\cite{Ma:mec}\cite{Ma:varen} (which included specific
proposals\cite{CWSSW:ten} in a general and systematic $R$-matrix form), this
programme has been stuck in recent years due to the following `dilaton
problem': when one tries to $q$-deform the $q$-Poincar\'e group one finds quite
generally\cite{Ma:poi} that, for the types of deformations of interest, one
must
include a scale generator as well. The purely Poincar\'e sector in this family
does not close as a Hopf algebra. A related problem is that even after this
extension, it does not seem possible to obtain a Hopf $*$-algebra (i.e. to
define complex conjugation of the generators) in any conventional sense.
Moreover, the extended $q$-Poincar\'e  Hopf algebra is not in general
quasitriangular, i.e. not a strict quantum group in the sense of
Drinfeld\cite{Dri}.

Here we solve these three problems as follows. First, we propose to embrace the
dilaton generator and $q$-deform the entire conformal algebra. This algebra is
isomorphic to the standard Drinfeld-Jimbo deformation $U_q(so(4,2))$ or a
cocycle twisting of it, but obtained now as an example of a new $R$-matrix
construction. We use the categorical double-bosonisation theory developed
recently in \cite{Ma:con}. From a physical point of view, this focus on the
conformal algebra limits our
theories at first to massless ones. This is not, however, an unreasonable
starting point, especially if we are interested in fundamental theories  where
observed particles are essentially massless compared to the Planck mass scale.
Secondly, we show that the natural $*$-structure on the $q$-conformal algebra
for real $q$ makes it into a {\em quasi}-$*$ Hopf algebra in the sense recently
introduced in \cite{Ma:qsta} for the Poincar\'e case. There is also a
quasitriangular structure. This work provides a first step towards the
development of a $q$-twistor theory based on the properties of $q$-Minkowski
space and the $q$-conformal algebra.

For quantum groups and braided groups, we adopt the conventions and notation
in \cite{Ma:book} and \cite{Ma:introp}. Briefly, a quantum group $H$ has a
coproduct $\Delta:H\to H\tens H$ which is a homomorphism to the usual tensor
product. By contrast, a braided group $B$ has a map $\und\Delta:B\to B\und\tens
B$ which is a homomorphism to a braided or non-commuting tensor product. This
concept and many examples have been introduced by the author\cite{Ma:exa}. In
physical terms, the elements of $B$ enjoy braid statistics. When discussing
quasitriangular structures, we will require formal power series in a
deformation
parameter in the usual way. All other constructions are algebraic.

Although we emphasise the construction in Proposition~2.1 below as a
$q$-conformal  algebra, we can also choose $R$ from other standard families
such as $su_n,sp_n$. Then the construction takes us up one in the family, i.e.
it allows the construction of quantum groups by induction\cite{Ma:con}. Or we
can choose non-standard $R$ and obtain altogether new quantum groups. Also,  it
has been pointed out to us that there are some superficial similarities between
some of the relations in Proposition~2.1 and some of the relations
independently proposed for a Hopf algebra construction in \cite{Vla:som}.  The
two constructions are, however,  not at all the same.

\subsection*{Acknowledgements} The contents of this paper were presented at the
minisemester on quantum groups and quantum spaces, Warsawa, Poland, November
and December 1995; I thank the organisers for their hospitality.

\section{$q$-Conformal algebra}

In this preliminary section we define a quasitriangular Hopf algebra which we
call the conformal algebra associated to regular $R$-matrix data. When $R$ is
the $so_n$ $R$-matrix, one obtains $U_q(so_{n+2})$. The construction is a
specific example of a left-module version of a recent abstract construction in
\cite{Ma:con}. This is explained further in the Appendix; here we describe only
the resulting algebra. In fact, the construction of the dilaton-extended
$q$-Poincar\'e algebra in $R$-matrix form has already been obtained (by the
author) in \cite{Ma:poi} and \cite{Ma:qsta}. We use exactly the results and
conventions developed  for this, extending it now by the special conformal
transformations.

Thus, let ${\bf R},{\bf R'}\in M_n\tens M_n$ be invertible matrices obeying the
conditions in \cite{Ma:poi}  for the construction of a linear (momentum)
braided group $V({\bf R'},{\bf R})$ and its associated extended $q$-Poincar\'e
algebra, which we denote $P({\bf R'},{\bf R})$.  The momentum sector has
generators $p^i$ and relations and braid statistics\cite{Ma:poi}
\eqn{mom}{ \vecp_1\vecp_2={\bf R'}\vecp_2\vecp_1,\quad \vecp_1'\vecp_2={\bf
R}\vecp_2\vecp_1'}
while $\vecp'$ denotes the generators of the second copy. We use a standard
compact notation where suffices denote tensor contraction positions. The
braided coproduct $\und\Delta \vecp=\vecp\tens 1+1\tens\vecp\equiv\vecp+\vecp'$
is called   a {\em braided coaddition}. The braided antipode is $\und
S\vecp=-\vecp$.
The rotation sector has matrix generators $l^\pm{}^i{}_j$ and the usual
relations\cite{FRT:lie}
\eqn{lor}{ \vecl^\pm_1\vecl^\pm_2{\bf R}={\bf R}\vecl^\pm_2\vecl^\pm_1,\quad
\vecl^+_1\vecl^-_2{\bf R}={\bf R}\vecl^-_2\vecl^+_1.}
{\em We also require} other relations beyond these quadratic ones, such that
the $\vecl^\pm$ with matrix coproduct $\Delta\vecl^\pm=\vecl^\pm\tens\vecl^\pm$
define a quasitriangular Hopf algebra under which (\ref{mom}) and (\ref{spcon})
remain covariant. For the $so_n$ R-matrix the $\vecl^\pm$ generate $U_q(so_n)$
in the FRT form \cite{FRT:lie}. Our construction is not limited to this case,
however. The dilatation sector is an additional generator $\dila$ with
coproduct $\Delta\dila=\dila\tens\dila$. The $q$-Poincar\'e algebra is
generated by these subalgebras, with the cross relations\cite{Ma:poi}
\eqn{poinrel}{ \vecl^+_1\vecp_2=\lambda^{-1} {\bf
R}_{21}^{-1}\vecp_2\vecl^+_1,\quad \vecl^-_1\vecp_2=\lambda{\bf
R}\vecp_2\vecl^-_1, \quad \dila \vecp=\lambda^{-1}\vecp \dila, \quad
[\vecl^\pm,\dila]=0,}
where $\lambda$ is the {\em quantum group normalisation constant}\cite{Ma:poi}
appearing in the fundamental representation $\rho(\vecl^+)=\lambda{\bf R}$,
$\rho(\vecl^-)=\lambda^{-1}{\bf R}^{-1}_{21}$ of the rotation sector. The
braided coproduct and antipode of the $\vecp$ generators do not appear directly
in the extended $q$-Poincar\'e algebra (which is an ordinary Hopf algebra), but
in the bosonised form
\eqn{poincop}{ \Delta \vecp=\vecp\tens 1+\dila\vecl^-\tens \vecp,\quad
S\vecp=-(S\vecl^-)\dila^{-1}\vecp.}

To this construction, we add now the special conformal transformations  as the
linear braided group $\Vhaj({\bf R'},{\bf R})^{\rm op}$ with generators $c_i$
and relations and braid statistics
\eqn{spcon}{\vecc_2\vecc_1=\vecc_1\vecc_2{\bf R'},\quad
\vecc_2'\vecc_1=\vecc_1\vecc_2'{\bf R}.}
There is a linear braided coproduct $\und\Delta\vecc=\vecc+\vecc'$.

\begin{propos} The extended $q$-Poincar\'e algebra (\ref{mom})--(\ref{poincop})
 has a further extension by generators $c_i$ obeying
\cmath{ \vecc_2\vecc_1=\vecc_1\vecc_2{\bf R'},\quad \vecl^+_1\vecc_2=\lambda
\vecc_2  \vecl^+_1{\bf R}_{21},\quad \vecl^-_1\vecc_2=\lambda^{-1}
\vecc_2\vecl^-_1{\bf R}^{-1}\\
 \quad \dila \vecc=\lambda\vecc \dila,\quad
[\vecp,\vecc]={\vecl^+\dila^{-1}-\vecl^-\dila\over q-q^{-1}},\quad \Delta
\vecc=\vecc\tens \vecl^+\dila^{-1}+1\tens\vecc,\quad S\vecc=-\vecc\dila
S\vecl^+}
where it is assumed that ${\bf R}$ depends on a parameter $q$ such that ${\bf
R}_{21}{\bf R}=\id+O(q-q^{-1})$. We call this the $q$-conformal algebra $C({\bf
R'},{\bf R})$ associated to our $R$-matrix data.
\end{propos}
\proof An abstract derivation is in the Appendix, but a direct proof is also
possible. Indeed, it is clear that $\vecc,\vecl^\pm,\dila$ generate a
`conjugate' Hopf algebra to the extended $q$-Poincar\'e one: their relations
are analogous under a symmetry. Hence it suffices to verify that the coproduct
is compatible with the cross relations. Thus,
\align{{}[\Delta p^i,\Delta c_j]&=&[p^i,c_a]\tens l^+{}^a{}_j\dila^{-1}+\dila
l^-{}^i{}_a\tens[p^a,c_j]+\dila l^-{}^i{}_a c_b\tens p^al^+{}^b{}_j\dila^{-1}-
 c_b\dila l^-{}^i{}_a\tens l^+{}^b{}_j\dila^{-1} p^a\\
&=&[p^i,c_a]\tens l^+{}^a{}_j\dila^{-1}+\dila l^-{}^i{}_a\tens[p^a,c_j]+\dila
c_dl^-{}^i{}_c {\bf R}^{-1}{}^c{}_a{}^d{}_b\tens p^al^+{}^b{}_j\dila^{-1}-
 c_b\dila l^-{}^i{}_a\tens l^+{}^b{}_j\dila^{-1} p^a\\
&=&{l^+{}^i{}_a\dila^{-1}\tens l^+{}^a{}_j\dila^{-1}-l^-{}^i{}_a\dila\tens
l^-{}^a{}_j\dila\over
q-q^{-1}}=\Delta{l^+{}^i{}_j\dila^{-1}-l^-{}^i{}_j\dila\over
q-q^{-1}}=\Delta[p^i,c_j]}
as required. We used the stated $\vecc,\vecl^-$ and $\vecp,\vecl^+$ relations
for the second and third equalities, as well as the $[\vecp,\vecc]$ relations
for the latter.

Note that the $q-q^{-1}$ factor ensures that our algebra has a reasonable limit
as $q\to 1$ but is not needed for the Hopf algebra structure itself (any factor
will do for this). \endproof

In a setting where $q=e^{t\over 2}$, there is typically a quasitriangular
structure $\CR_L$ for the Lorentz/rotation sector as a formal power series in
$t$. In this setting:

\begin{propos}  The $q$-conformal Hopf algebra above is   quasitriangular, with
\eqn{confqua}{\CR=\CR_L \lambda^{-\xi\tens\xi} \exp(\vecc|\vecp)^{-1}}
where $\CR_L$ is the $q$-Lorentz quantum group quasitriangular structure,
$\dila=\lambda^\xi$ and $\exp(\vecc|\vecp)\in \Vhaj\tens V$ is the {\em braided
exponential} or canonical element associated with the braided group duality
pairing between $V$ and $\Vhaj$ as linear braided groups.
\end{propos}
\proof This follows from the general construction\cite{Ma:con} underlying the
above proposition; see the Appendix. To verify it directly, one may use  the
bicharacter property of the braided exponentials under the braided
coproduct\cite{Ma:fre}, with the corresponding properties with respect to the
bosonised coproduct.  \endproof

The braided exponential $\exp(\vecc|\vecp)$ here is a $\vecp$-eigenfunction or
plane wave in the copy of $q$-spacetime generated by the $c_i$ (and likewise a
plane wave in the copy generated by $p^i$).

\section{Quasi $*$-structure on the conformal generators}

So far, we have considered the complexified picture. We now consider
$*$-structures on our algebras. The specification of the $*$-structure in the
momentum sector determines which linear combinations are `real' in the sense of
being invariant under $*$. This determines which representations are unitary
(such elements should be Hermitian) and also determines, when there is a
quantum metric, whether it is of Euclidean, Minkowski or other type  according
to the form of its restriction to the such elements. $*$-Structures for the
extended $q$-Poincar\'e group have already been analysed in \cite{Ma:qsta}, and
we extend this now to the $q$-conformal case. We recall from \cite{Ma:qsta}
that one needs the notion of a quasi-$*$ Hopf algebra $\CH$. This is a Hopf
algebra over $\C$ which is a $*$-algebra, and an invertible element $\CS\in
\CH\tens\CH$ such that
\eqn{qsta}{ (*\tens *)\circ\Delta\circ *=\CS^{-1}(\tau\circ\Delta\ )\CS,\quad
(\Delta\tens\id)\CS=\CS_{13}\CS_{23},\quad
(\id\tens\Delta)\CS=\CS_{13}\CS_{12},\quad \CS^{*\tens *}=\CS_{21},}
where $\tau$ denotes transposition. One can show that $\CS$ obeys the QYBE (but
we do not denote it here by $\CR$, to avoid confusion with the quasitriangular
structure also present below). To have such a structure in our R-matrix setting
we suppose that the R-matrix in the preceding section is of one of the two real
types in \cite{Ma:star}. We also suppose a quantum metric $\eta$ compatible
with ${\bf R}$ (see \cite{Ma:varen}) and of corresponding reality type:
\eqn{Rtpye}{ \overline{{\bf R}^i{}_j{}^k{}_l}=\cases{{\bf R}^l{}_k{}^j{}_i&
Real\ Type\ I\cr
{\bf R}^{\bar j}{}_{\bar l}{}^{\bar i}{}_{\bar k}&Real\ Type\ II},\quad
\overline{\eta_{ij}}=\cases{\eta^{ji}&Real\ Type\ I\cr \eta_{\bar j\bar
i}&Real\ Type\ II},}
where $\eta^{ij}$ is the transposed inverse of $\eta_{ij}$ and $\bar i$ is an
involution on the indices assumed in the type II case. We assume
$\lambda^*=\lambda$ as well. These reality conditions hold for the standard
choices  ${\bf R}$, when their parameter $q$ is real.

The extended $q$-Poincar\`e algebra has the quasi-$*$-structure\cite{Ma:qsta}
\eqn{pstarI}{ p^i{}^*=\cases{\eta_{ia}p^a&Real Type I\cr p^{\bar i}&Real\ Type\
II},\quad l^\pm{}^i{}_j{}^*=\cases{\eta_{ib}l^\mp{}^b{}_a\eta^{ja}&Real\ Type\
I\cr l^\mp{}^{\bar i}{}_{\bar j}&Real\ Type\ I},\quad \dila^*=\dila^{-1},\quad
\CS=\CR_{L}\lambda^{-\xi\tens\xi}}
where $\CR_{L}\lambda^{-\xi\tens\xi}$ is the dilaton-extended $q$-Lorentz
quasitriangular structure and $\dila=\lambda^\xi$. Note that  the real type II
case used in \cite{Ma:qsta} was chosen such that on the Lorentz generators in
function-algebra form it appears as $t^i{}_j{}^*=t^{\bar i}{}_{\bar j}$, which
corresponds to $l^\pm{}^i{}_j{}^*=S^2l^\mp{}^{\bar i}{}_{\bar j}$. We can
equivalently put the $S^2$ automorphism on the function algebra side as
$t^i{}_j{}^*=S^2t^{\bar i}{}_{\bar j}$, as we prefer now.

\begin{propos} The quasi-$*$ structure (\ref{pstarI}) extends to one on the
$q$-conformal algebra in Proposition~2.1, with
$c_i{}^*=\cases{c_a\eta^{ia}&Real\ Type\ I\cr c_{\bar i} &Real\ Type\ II}$.
Moreover,
\[ (*\tens *)\circ\Delta\circ *=\exp(\vecc|\vecp)^{-1}(\Delta\
)\exp(\vecc|\vecp)\]
holds for the coproduct on any element of the $q$-conformal algebra.
\end{propos}
\proof The proof of compatibility of this $*$ with the $\vecc,\vecl^\pm$
relations is similar to that for $\vecp,\vecl^\pm$. Explicitly, in the type I
case,
\align{(\lambda  c_b l^+{}^i{}_a {\bf R}^b{}_k{}^a{}_j)^*\equad &&=\lambda {\bf
R}^j{}_a{}^k{}_b\eta^{bd}\eta_{ic}l^-{}^c{}_e\eta^{ae}c_d\\
&&=\lambda {\bf R}^j{}_a{}^k{}_b\eta^{bd}\eta_{ic}\eta^{ae}\lambda^{-1} c_g
l^-{}^c{}_f {\bf
R}^{-1}{}^f{}_e{}^g{}_d=c_a\eta^{ka}\eta_{ic}l^-{}^c{}_d\eta^{jd}=(l^+{}^i{}_j
c_k)^*}
using invariance of $R$ under conjugation by $\eta\tens\eta$. In addition, we
have
\[
[p^i,c_j]^*=[c_j{}^*,p^{i*}]=[c_a\eta^{ja},p^b\eta_{ib}]=-\eta^{ja}\eta_{ib}({
l^+{}^b{}_a\dila^{-1}-l^-{}^b{}_a\dila\over
q-q^{-1}})=({l^+{}^i{}_j\dila^{-1}-l^-{}^i{}_j\dila\over q-q^{-1}})^*\]
as required. Hence we have a $*$-algebra in this case. In the type II case, the
calculation is
\[(\lambda  c_b l^+{}^i{}_a {\bf R}^b{}_k{}^a{}_j)^*=\lambda l^-{}^{\bar
i}{}_{\bar a} c_{\bar b} {\bf R}^{\bar a}{}_{\bar j}{}^{\bar b}{}_{\bar
k}=c_{\bar k}l^-{}^{\bar i}{}_{\bar j}=(l^+{}^i{}_j c_k)^*\]
\[ [p^i,c_j]^*=[c_{\bar j},p^{\bar i}]=-({l^+{}^{\bar i}{}_{\bar
j}\dila^{-1}-l^-{}^{\bar i}{}_{\bar j}\dila\over
q-q^{-1}})=({l^+{}^i{}_j\dila^{-1}-l^-{}^i{}_j\dila\over q-q^{-1}})^*.\]
In either case, the sub-Hopf algebra generated by $\dila,\vecc,\vecl^\pm$ forms
a quasi-$*$ Hopf algebra
with the same cocycle $\CS=\CR_L\lambda^{-\xi\tens\xi}$, by analogous arguments
to the proof for the extended $q$-Poincar\'e algebra in \cite{Ma:qsta}.
Combining (\ref{qsta}) with Proposition~2.2 gives the form of $(*\tens
*)\circ\Delta\circ *$ stated. \endproof

Although the $q$-conformal algebra with the above $*$-operation is not a Hopf
$*$-algebra in the usual sense, we see that $\Delta$ fails to be a
$*$-algebra map only up to conjugation by the plane-wave $\exp(\vecc|\vecp)$.
More precisely,  every quasi-$*$ Hopf algebra the conjugate coproduct
$\bar\Delta=(*\tens *)\circ\Delta\circ *$  also provides a quasi-$*$ Hopf
algebra structure, in general different from
$\Delta$. In our case, this comes out as
\eqn{barDelta}{\bar\Delta \vecc=\vecc\tens\vecl^-\dila+1\tens\vecc}
for either the real type I or type II $*$-structures above (similarly for
$\bar\Delta\vecp$ in \cite{Ma:qsta}); Proposition~3.1 tells us that this
$\bar\Delta$  and the coproduct $\Delta$ in Proposition~2.1 are conjugate by
$\exp(\vecc|\vecp)$.

\section{Spinorial formulation}

An important class of examples  of our data
 ${\bf R'},{\bf R}$  is provided by a `spinorial' construction starting from a
smaller Yang-Baxter matrix $R\in M_s\tens M_s$, where $n=s^2$. We require this
to be $q$-Hecke in the sense $(PR-q)(PR+q^{-1})=0$, where $P$ is the
permutation matrix. The extended $q$-Poincar\'e algebra in this setting has
been given in \cite{Ma:qsta}, while the momentum sector or $q$-spacetime itself
is from \cite{Ma:euc}\cite{Ma:exa}\cite{Mey:new} and  reviewed in
\cite{Ma:varen} or \cite[Ch.~10]{Ma:book}. We include now the $q$-conformal
algebra in this spinorial approach. In fact, the construction has two versions
which are strictly `gauge equivalent' in a certain algebraic sense.  These are
the `Euclidean' and `Minkowski' gauges of the same construction introduced in
\cite{Ma:euc} and \cite{Ma:exa}\cite{Mey:new} respectively.

The Euclidean gauge construction is \cite{Ma:euc}
\eqn{eucR}{ {\bf R'}^i{}_j{}^k{}_l=R^{-1}{}^{l_0}{}_{k_0}{}^{j_0}{}_{i_0}
R^{i_1}{}_{j_1}{}^{k_1}{}_{l_1},\quad {\bf
R}^i{}_j{}^k{}_l=R{}^{j_0}{}_{i_0}{}^{l_0}{}_{k_0}
R^{i_1}{}_{j_1}{}^{k_1}{}_{l_1}}
and is equivalent for $R$ the standard $su_2$ R-matrix to taking for ${\bf R}$
the standard $so_4$ R-matrix. Of course, the construction is more general and
can be used just as well to define non-standard spacetimes by taking other
non-standard $R$.  We write $i=i_0i_1$, $j=j_0j_1$ etc as multi-indices.

We also write $p^i=p^{i_1}{}_{i_0}$. Then the relations (\ref{mom}) in the
momentum sector (which will also be the relations of $q$-spacetime) become
\cite{Ma:euc}
\eqn{sseuc}{ \quad R_{21}\vecp_1\vecp_2=\vecp_2\vecp_1R.}
More non-trivially, we replace the vectorial $q$-Lorentz algebra  generated by
$l^\pm{}^i{}_j$ by a spinorial version generated by two sets of generators
$l^\pm{}^{i_0}{}_{j_0}$ and $m^\pm{}^{i_1}{}_{j_1}$ obeying relations like
(\ref{lor}) with respect to $R$. For $R$ the $su_2$ R-matrix, the momentum and
spacetime sectors are isomorphic to the quantum matrices $\bar M_q(2)$, and the
Lorentz/rotation sector is $U_q(su_2)\tens U_q(su_2)$. The natural
$*$-structure in this gauge is the unitary type one which corresponds to
$SU_q(2)$ as a $q$-deformed 3-sphere in $M_q(2)$. The dilaton sector is
generated by $\dila$ as before, commuting with $\vecl^\pm,\vecm^\pm$. The
cross-relations between these various sectors and the coproducts are obtained
in \cite{Ma:qsta}. In the present (slightly different) conventions they come
out as:
\ceqn{eucpoi}{\vecp_1\vecl^+_2=\lambda^{-\h}R^{-1}\vecl^+_2\vecp_1,\quad
\vecp_1\vecl^-_2=\lambda^{\h}R_{21}\vecl^+_2\vecp_1,\quad
\vecp_1\vecm^+_2=\lambda^{\h}R\vecm^+_2\vecp_1,\quad
\vecp_1\vecm^-_2=\lambda^{-\h}R^{-1}_{21}\vecm^-_2\vecp_1,\\
\dila\vecp=\lambda^{-1}\vecp\dila,\quad \Delta\vecp=\vecp\tens 1+\dila
S^{-1}(S\vecm^-(\ )\vecl^-)\tens \vecp,\quad \eps(\vecp)=0}
where the space is for the matrix indices of $\vecp$ to be inserted.

To this spinorial extended $q$-Poincar\'e algebra, we add the special conformal
transformations $c_i=c^{i_0}{}_{i_1}$. Note that the assignment is transposed
relative to the assignment for $p^i$.

\begin{propos} In the Euclidean gauge, the spinorial extended $q$-Poincar\'e
algebra in \cite{Ma:qsta} has a further extension by a matrix of generators
$\vecc$ obeying
\cmath{ R\vecc_1\vecc_2=\vecc_2\vecc_1R_{21},\quad
\vecl^+_1\vecc_2=\lambda^{-\h}\vecc_2\vecl^+_1R^{-1},\quad
\vecl^-_1\vecc_2=\lambda^{\h}\vecc_2\vecl^-_1R_{21},\\
\vecm^+_1\vecc_2=\lambda^{\h}\vecc_2\vecm^+_1R,\quad
\vecm^-_1\vecc_2=\lambda^{-\h}\vecc_2\vecm^-_1R^{-1}_{21},\quad
\dila\vecc=\lambda\vecc\dila,\\
{}[p^{i_1}{}_{i_0},c^{j_0}{}_{j_1}]={\dila^{-1}(S^{-1}l^+{}^{j_0}{}_{i_0})
m^+{}^{i_1}{}_{j_1}-\dila (S^{-1}l^-{}^{j_0}{}_{i_0})m^-{}^{i_1}{}_{j_1}
\over q-q^{-1}},\quad
\Delta\vecc=\vecc\tens \dila^{-1}(S^{-1}\vecl^+)(\ )\vecm^++1\tens \vecc,\quad
\eps(\vecc)=0}
and forming a quasitriangular Hopf algebra. This is the spinorial
$q$-conformal algebra in the Euclidean gauge.
\end{propos}
\proof The $\vecc,\vecl^\pm,\vecm^\pm,\dila$ relations are obtained along the
same lines as in \cite{Ma:qsta} via double-bosonisation. They are consistent
with Proposition~2.1 using (\ref{eucR}) and the ansatz
$l^\pm{}^i{}_j=(S^{-1}l^\pm{}^{j_0}{}_{i_0})m^\pm{}^{i_1}{}_{j_1}$. The
$\vecc,\vecp$ relations and the coproduct follow at once
from this form of $l^\pm{}^i{}_j$. We note that if we use the (slightly
different) identification $\eta_{ia}p^a=p^{i_0}{}_{i_1}$ used in
\cite{Ma:qsta}, and the expression $\eta_{ij}=\eps^{i_0j_0}\eps_{i_1j_1}$ in
terms of the spinor metric associated to $R$, then the $[p^i,c_j]$ relations
come out as
\eqn{pcold}{[p^{i_0}{}_{i_1},c^{j_0}{}_{j_1}]={\dila^{-1}l^+{}^{i_0}{}_{a_0}
\eps^{a_0j_0} \eps_{a_1i_1} m^+{}^{a_1}{}_{j_1}-\dila l^-{}^{i_0}{}_{a_0}
\eps^{a_0j_0} \eps_{a_1i_1} m^-{}^{a_1}{}_{j_1}\over q-q^{-1}}.}
The spinor metric also converts the $*$-structure in Section~3
into a matrix form. \endproof

The Minkowski gauge for the same construction is\cite{Ma:exa}\cite{Mey:new}
\eqn{minkR}{{\bf R'}^i{}_j{}^k{}_l=R^{-1}{}^{d}{}_{k_0}{}^{j_0}{}_{a}
R^{k_1}{}_{b}{}^{a}{}_{i_0}R^{i_1}{}_c{}^b{}_{l_1} {\widetilde
R}^c{}_{j_1}{}^{l_0}{}_d,\quad {\bf R}{}^i{}_j{}^k{}_l=R^{j_0}{}_a{}^d{}_{k_0}
R^{k_1}{}_b{}^a{}_{i_0}
R^{i_1}{}_c{}^b{}_{l_1} {\widetilde R}^c{}_{j_1}{}^{l_0}{}_d.}
The momentum or spacetime sector in this case has the braided matrix
relations\cite{Ma:exa}
\eqn{ssmink}{R_{21}\vecp_1 R\vecp_2=\vecp_2 R_{21}\vecp_1 R}
where $\eta_{ia}p^a=p^{i_0}{}_{i_1}$, and yields the braided matrices $BM_q(2)$
for the standard $su_2$ R-matrix. The natural spacetime $*$-structure in this
case is a Hermitian one, justifying the name for this gauge. (The unit sphere
here is actually isomorphic to $U_q(su_2)$ as a $*$-algebra when $q\ne 1$). The
Lorentz sector in this standard case is $U_q(su_2)\codcross U_q(su_2)$ (with a
more complicated coproduct than in the Euclidean gauge).

The Euclidean gauge for $q$-spacetime was introduced in \cite{Ma:euc} precisely
as gauge equivalent to the  Minkowski gauge (which was found first). At the
Lorentz algebra level the gauges are related by twisting by a quantum cocycle
(see \cite[Sec. 4]{Ma:poi}, in a dual form). This was extended to the level of
the extended $q$-Poincar\'e algebra in \cite{Ma:qsta}, using the same cocycle
viewed in the bigger algebra. The cocycle is $\chi=\CR_{23}^{-1}$ where $\CR$
is the quasitriangular structure of $U_q(su_2)$ in the standard example.

\begin{propos} The same quantum cocycle $\chi$ viewed in the spinorial form of
the $q$-conformal algebra twists its structure from the Euclidean to the
Minkowski gauge.
\end{propos}
\proof This is true for the sub-Hopf algebra generated by
$\vecc,\vecl^\pm,\dila$ by analogous arguments to those for the extended
$q$-Poincar\'e algebra. Since the coproduct is entirely defined by its
restriction to either of these two sub-Hopf algebras, we conclude the same
twisting result for the entire $q$-conformal algebra. \endproof

In view of this, we will not give the structure in detail in the Minkowski
gauge: the structure of the spinorial form of the extended $q$-Poincar\'e
algebra is given   in \cite{Ma:qsta}. To this, we add the special conformal
transformations in the form
  $\bar R_{21}\vecc_1\bar R\vecc_2=\vecc_2\bar R_{21}\vecc_1 \bar R$, where
$c_i=c^{i_1}{}_{i_0}$ and $\bar R^i{}_j{}^k{}_l=R^l{}_k{}^j{}_i$. The
cross relations with $\vecl^\pm$ are similar to those between $\vecp$ and
$\vecl^\pm$ in \cite{Ma:qsta}.

\section{Conformal transformations of spacetime}

So far, we have called our quasitriangular Hopf algebra $C({\bf R'},{\bf R})$
the $q$-conformal one because of its structural form, which is analogous to
that of the conformal Lie algebra. We are now ready to justify the terminology
in geometrical terms, i.e. by its action on $q$-spacetime. For the latter, we
take the linear braided group $\Vhaj({\bf R'},{\bf R})$ with generators $x_i$
and relations and braid statistics
\eqn{space}{ \vecx_1\vecx_2=\vecx_2\vecx_1{\bf R'},\quad
\vecx'_1\vecx_2=\vecx_2\vecx'_1{\bf R}.}
There is a linear coproduct $\und\Delta\vecx=\vecx+\vecx'$ and a $*$-structure
$x_i^*=\cases{x_a\eta^{ia}&Real\ Type\ I\cr x_{\bar i}&Real\ Type II}$,
which we take of the same form as for $\vecc$ in Section~3.

{}From the theory of braided groups, it is known\cite{Ma:poi} that the extended
$q$-Poincar\'e algebra acts covariantly on $q$-spacetime by $q$-rotations (via
the fundamental representation defined by ${\bf R}$) and
braided-differentiation
for the momentum sector\cite{Ma:poi}\cite{Ma:qsta}
\eqn{actpoi}{\vecl^+_1\la \vecx_2=\vecx_2\lambda {\bf R}_{21},\quad
\vecl^-_1\la \vecx_2=\vecx_2\lambda^{-1}{\bf R}^{-1},\quad p^i\la
x_j=-\delta^i{}_j,\quad \dila\la x_i=\lambda x_i.}
To this, we add:

\begin{propos} The $q$-conformal algebra in Proposition~2.1 acts covariantly on
$q$-spacetime by (\ref{actpoi}) and
\[ \vecc_2\la \vecx_1={\vecx_1\vecx_2-\vecx_2\vecx_1{\bf R}\over q-q^{-1}}.\]
\end{propos}
\proof This follows from general theory in \cite{Ma:con}; the required action
of $\vecc$ is derived in the Appendix. The direct proof that the
$\vecc,\vecl^\pm,\dila$ relations are represented is similar to that for
$\vecp,\vecl^\pm,\dila$. For the $\vecp,\vecc$ relations we can check it easily
at lowest order, as $(q-q^{-1})[p^i,c_j]\la x_k=p^i\la (x_kx_j- x_bx_a{\bf
R}^a{}_k{}^b{}_j) + c_j\la\delta^i{}_k=-\delta^i{}_k x_j - x_a {\bf
R}^{-1}{}^i{}_j{}^a{}_k+\delta^i{}_b x_a {\bf R}^a{}_k{}^b{}_j+x_c{\bf
R}^{-1}{}^i{}_a{}^c{}_b {\bf R}^a{}_k{}^b{}_j=(\dila^{-1}l^+{}^i{}_j - \dila
l^-{}^i{}_j)\la x_k$, where the outer two terms cancelled. We used the action
of $p^i$ on products $x_jx_k$ via the braided-Leibniz rule with ${\bf
R}_{21}^{-1}$ \cite{Ma:qsta}. One can proceed similarly for the higher order
case, using the action of $c_i$ on products obtained below.
\endproof

Note that both the action of $p^i$ and $c_i$  extend to products via a
braided-Leibniz rule, because they originate as braided module algebra
structures (this is equivalent to the statement that the actions form
a module-algebra structure with respect to the Hopf algebra coproducts.) In the
case of $p^i$, the action on a general monomial comes out in terms of the
braided-integer matrices with respect to ${\bf R}^{-1}_{21}$ (see
\cite{Ma:qsta}).
For the $c_i$ we have:

\begin{lemma} The action of $c_i$ on a general product is
\[
\vecc_n\la\vecx_1\vecx_2\cdots\vecx_{n-1}=\vecx_1\vecx_2\cdots\vecx_n({1-(P{\bf
R})_{12}(P{\bf R})_{23}\cdots (P{\bf R})_{n-1n}\over q-q^{-1}})\]
where $P$ is the permutation matrix.
\end{lemma}
\proof We first  compute the braided-Leibniz rule for $c_i$. As explained in
the Appendix, its natural form is as a right-handed (braided) derivation
$c_i\la =\delta_i$ acting from the right. Then
\[ (ab)\delta_i=a(b\delta_i)+a\Psi(b\tens\delta_i)\]
where the braiding is the braiding for the covector braided group $\Vhaj({\bf
R'},{\bf R})$, i.e. defined by ${\bf R}$. Hence
\[ (\vecx_1\cdots\vecx_{n-1})\delta_n=\vecx_1\cdots \vecx_n({1-(P{\bf R})_{n-1
n}\over q-q^{-1}})+(\vecx_1\cdots\vecx_{n-2})\delta_{n-1}\vecx_n (P{\bf
R})_{n-1
n}.\]
The result then follows by induction. \endproof

Another way to describe the action is in terms of the algebra structure
of the corresponding semidirect product  of spacetime crossed by the
$q$-conformal group. The cross relations between the extended $q$-Poincar\'e
algebra and spacetime is \cite{Ma:qsta}
\eqn{poincx}{ \vecl^+_1\vecx_2=\vecx_2 \lambda {\bf R}_{21} \vecl^+_1,\quad
\vecl^-_1\vecx_2=\vecx_2 \lambda^{-1} {\bf R}^{-1} \vecl^-_1,\quad \vecx_2{\bf
R}^{-1}\vecp_1-\vecp_1\vecx_2=\id,\quad  \dila\vecx=\lambda \vecx\dila.}
The $\vecx,\vecp$ relations are the `braided Heisenberg algebra' in the present
conventions. To this we now add:

\begin{propos} The  $q$-conformal group acting as above and $\vecx$ acting
by left multiplication on $q$-spacetime form a representation of the algebra
$\Vhaj({\bf R'},{\bf R})\lcross C({\bf R'},{\bf R})$ with the additional
$\vecc,\vecx$ cross relations
\[ {}[\vecc_1+{\vecx_1\vecl^+_1\dila^{-1}\over q-q^{-1}},\vecx_2]=0.\]
\end{propos}
\proof We make a left handed semidirect product using the coproduct in
Proposition~2.1, the action of $c_i$ above and the already-known
cross-relations (\ref{poincx}). Thus
\align{c_ix_j\equad && =(c_i\o\la x_j)c_i\t=(c_a\la
x_j)\vecl^+{}^a{}_i\dila^{-1}+x_j c_i\\
&&={x_jx_a - x_d x_c R^c{}_j{}^d{}_a\over
q-q^{-1}}l^+{}^a{}_i\dila^{-1}+x_jc_i={x_jx_a l^+{}^a{}_i\dila^{-1}-x_a
l^+{}^a{}_i\dila^{-1} x_j \over q-q^{-1}}+x_jc_i,}
as stated. Because the action in Proposition~5.1 is covariant ($q$-spacetime
forms
a module algebra under it), we know from the general theory of Hopf algebra
cross products that these relations define an associative algebra structure on
the tensor product vector space, and that the action on $q$-spacetime extends
to it with $x_i$ acting by left-multiplication. \endproof

We can also use the spinorial form of the $q$-conformal algebra. The action of
the spinorial form of the extended $q$-Poincar\'e algebra is given in
\cite{Ma:qsta}. To this, we add:

\begin{propos} The spinorial form of the $q$-conformal algebra in the Euclidean
gauge acts as in \cite{Ma:qsta} and
\[ \vecc_2\la\vecx_1=-\vecx_1\vecx_2 PR.\]
\end{propos}
\proof We use the form of ${\bf R}$ in (\ref{eucR}) in Proposition~5.1 and
$R=R_{21}^{-1}+(q-q^{-1})P$ from the $q$-Hecke assumption in Section~4. Thus
\align{(q-q^{-1})c_j\la x_i\equad &&= x_i x_j-x_bx_a{\bf R}^a{}_i{}^b{}_j\\
&&=x^{i_0}{}_{i_1}x^{j_0}{}_{j_1}-x^{b_0}{}_{b_1} x^{a_0}{}_{a_1}
R^{-1}{}^{i_0}{}_{a_0}{}^{j_0}{}_{b_0} R^{a_1}{}_{i_1}{}^{b_1}{}_{j_1} -
(q-q)^{-1} x^{b_0}{}_{b_1} x^{a_0}{}_{a_1} \delta^{i_0}_{b_0}\delta^{j_0}_{a_0}
R^{a_1}{}_{i_1}{}^{b_1}{}_{j_1}.}
The first two terms then give zero due to the form of ${\bf R'}$ in
(\ref{eucR}) and the relations for the $x_i$. \endproof

We are now in position to compute this action for our standard
$q$-spacetime\cite{Ma:varen}.  The classical formula $c_j\la
x_i=\h\eta_{ij}\vecx\cdot\vecx-x_ix_j$ would be
\eqn{classcon}{ \pmatrix{\alpha&\beta\cr \gamma&\delta}\la\pmatrix{a&b\cr
c&d}=\pmatrix{
-a^2&-ba&-ac&-bc\cr
-ab&-b^2&-ad&-bd\cr
-ca&-da&-c^2&-dc\cr
-cb&-db&-cd&-d^2}}
where $\vecc=\pmatrix{\alpha&\beta\cr \gamma&\delta}$, $\vecx=\pmatrix{a&b\cr
c&d}$ and $\eta=\pmatrix{0&0&0&1\cr 0&0&-1&0\cr 0&-1&0&0\cr 1&0&0&0}$ is the
metric on complexified spacetime in these spinor coordinates (different linear
combinations are considered real spacetime coordinates in the Minkowski and
Euclidean cases).

\begin{example} For the standard $q$-spacetime in the Euclidean gauge, we have
\[ \pmatrix{\alpha&\beta\cr \gamma&\delta}\la\pmatrix{a&b\cr c&d}=\pmatrix{
-a^2&-q^2ab&-ca&-bc\cr
-ba&-b^2&-da&-db\cr
-ac&-ad-(q-q^{-1})bc&-c^2&-qdc\cr
-bc&-bd&-dc&-d^2}
\]
This is a $q$-deformation of the usual action of the special conformal
transformations  on spacetime.
\end{example}
\proof This is computed easily from Proposition~5.4 with
\[ R=\pmatrix{q&0&0&0\cr 0&1&q-q^{-1}&0\cr 0&0&1&0\cr 0&0&0&q}\]
which is the standard $su_2$ R-matrix in the $q$-Hecke normalisation.  The
relations between the non-commutative spinor
spacetime coordinates in this case are given explicitly in \cite{Ma:euc}.
The quantum group normalisation of the corresponding ${\bf R}$ is
$\lambda=q^{-1}$.
\endproof

This justifies our proposal for $C({\bf R'},{\bf R})$ as $q$-conformal group.
Note that the metric does not play any direct role in our definition of the
$q$-conformal group and its action on spacetime, i.e. our approach is a novel
one even when $q=1$. It is remarkable therefore that it coincides for our
standard example with the action (\ref{classcon}) defined through a metric.
The connection is quite general, however.

\begin{lemma} If $\eta$ is a quantum metric such that
$\vecx\cdot\vecx=x_ax_b\eta^{ba}$ is central, then it is preserved by the
$q$-conformal group up to scaling, in the sense
\[ c_i\la (\vecx\cdot\vecx)^m=({1-\lambda^{-2m}\over q-q^{-1}})x_i
(\vecx\cdot\vecx)^m,\quad {\rm i.e.} \quad c_i\la
f(\vecx\cdot\vecx)=({1-\lambda^{-2}\over q-q^{-1}})x_i \vecx\cdot\vecx
(\del_{\lambda^{-2}}f)(\vecx\cdot\vecx).\]
\end{lemma}
\proof We compute $(q-q^{-1})c_i\la (x_ax_b\eta^{ba})=x_ax_bx_i\eta^{ba}-x_c
x_d x_e {\bf R}^d{}_a{}^c{}_f {\bf R}^e{}_b{}^f{}_i\eta^{ba}=
x_ax_bx_i\eta^{ba}-x_c x_d x_e {\bf R}^d{}_a{}^c{}_f \eta^{eb}\lambda^{-2}{\bf
R}^{-1}{}^a{}_b{}^f{}_i=\vecx\cdot\vecx x_i-\lambda^{-2}x_i\vecx\cdot\vecx$,
using the covariance properties of the quantum metric. Under a further
condition on the quantum metric (true in the main examples, see
\cite{Ma:varen}) one knows that $\vecx\cdot\vecx$ is also central. This gives
the result for $m=1$. From the covariance properties of the quantum metric, we
likewise compute the braiding $\vecx'\cdot\vecx' x_i=x'_ax'_n x_i\eta^{ba}
=x'_a x_dx'_cR^c{}_b{}^d{}_i\eta^{ba}=x'_ax_dx'_c\lambda^{-2}
R^{-1}{}^a{}_b{}^d{}_i\eta^{cb}=x_fx'_ex'_c R^e{}_a{}^b{}_d\lambda^{-2}
R^{-1}{}^a{}_b{}^d{}_i\eta^{cb}=\lambda^{-2}x_i
\vecx'\cdot\vecx'$, i.e. $\Psi(\vecx\cdot\vecx\tens x_i)=\lambda^{-2}\vecx_i
\tens\vecx\cdot\vecx$. The $q$-Leibniz rule for the action of $c_i$ then 
implies $c_i\la(\vecx\cdot\vecx)^m=(\vecx\cdot\vecx)^{m-1}x_i({1-\lambda^{-2}
\over q-q^{-1}})\vecx\cdot\vecx+(c_i\la(\vecx\cdot\vecx)^{m-1})\vecx\cdot\vecx
\lambda^{-2}$, which provides the general result by induction. We alternatively 
write this in terms of a  $\lambda^{-2}$-deformed derivative defined in a 
usual way. \endproof

Thus, our $q$-conformal group and its action do not {\em a priori} involve a
metric,
but when there is one, it is preserved in some sense. Instead, the structures
and formulae which we normally associate with preservation up to scale of a
metric are obtained from the braided adjoint action. For example, we see that
the standard $q$-Gaussian $g_\eta$ in the setting of Lemma~5.6, which is a
$\lambda^{-2}$-exponential of $\vecx\cdot\vecx$, is preserved in the sense
\[ c_i\la g_\eta=-q^{-1}({1-\lambda^{-2}\over 1-q^{-4}})
x_i(\vecx\cdot\vecx)g_\eta,\]
in addition to its usual properties under the extended $q$-Poincar\'e algebra.

Finally, when a quasi-$*$ Hopf algebra acts covariantly on a $*$-algebra then
its conjugate quasi-$*$ Hopf algebra acts with a conjugate
action\cite{Ma:qsta}. In the case of the extended $q$-Poincar\'e algebra it was
 shown that the  action of $\vecl^\pm,\dila$ on $q$-spacetime is
self-conjugate, while the conjugate action of the $\vecp$ generators is by
braided-differentiation with ${\bf R}_{21}^{-1}$ replaced by ${\bf R}$.

\begin{propos} The conjugate action of the $q$-conformal algebra on
$q$-spacetime is
\[ \vecc_2\bar\la \vecx_1={\vecx_1\vecx_2-\vecx_2\vecx_1{\bf R}_{21}^{-1}\over
q-q^{-1}}.\]
Moreover, $c_i\la \und S(\ )=\und S(c_i\bar\la(\ ))$, where $\und S$ is the
braided antipode or parity operator on $q$-spacetime.
\end{propos}
\proof The abstract treatment for the conjugate action of the special conformal
generators is in the Appendix, from which one may compute the explicit form
stated.  In our R-matrix setting, a direct proof is as follows, using the
$*$-structures in Section~3 for either the real type I or type II cases. In the
type I case,
\align{c_i\bar\la x_j\equad && =(Sc_a \eta^{ia}\la x_b
\eta^{jb})^*=-(\eta^{ia}c_d \dila Sl^+{}^d{}_a\la x_b\eta^{jb})^*=-(\eta^{ia}
c_d\la x_e R^{-1}{}^e{}_b{}^d{}_a \eta^{jb})^*\\
&&=(({-x_e x_d + x_fx_g R^g{}_e{}^f{}_d\over
q-q^{-1}})\eta^{ia}R^{-1}{}^e{}_b{}^d{}_a\eta^{jb})^*={x_jx_i-x_a x_b
R^{-1}{}^a{}_i{}^b{}_j\over q-q^{-1}}}
using the usual covariance properties of the quantum metric. In the type II
case,
\[c_i\bar\la x_j=(Sc_{\bar i}\la x_{\bar j})^*=-(c_a \dila Sl^+{}^a{}_{\bar
i}\la x_{\bar j})^*=-(  c_a\la x_b R^{-1}{}^b{}_{\bar j}{}^a{}_{\bar
i})^*=(({-x_b x_a + x_dx_c R^c{}_b{}^d{}_a\over q-q^{-1}}) R^{-1}{}^b{}_{\bar
j}{}^a{}_{\bar i})^*\]
which likewise computes to the stated formula.

For the result that the action and conjugate action are intertwined by $\und
S$, we have on the generators $\und S(\vecc_2\bar
\la\vecx_1)=(q-q^{-1})^{-1}\und S(\vecx_1\vecx_2-\vecx_2\vecx_1{\bf
R}^{-1}_{21})=(q-q^{-1})^{-1}((-\vecx_2)(-\vecx_1) {\bf
R}-(-\vecx_1)(\vecx_2))=\vecc_2\la (-\vecx_1)$ using the
braided-antimultiplicativity of the braided antipode $\und S(\vecx)=-\vecx$.
\endproof

The same applies to the action on products of $q$-spacetime generators: we use
${\bf R}^{-1}_{21}$ in place of ${\bf R}$ in Lemma~5.2. Thus the $q$-conformal
group  exhibits the same novel phenomenon  demonstrated for the extended
$q$-Poincar\'e algebra in \cite{Ma:qsta} whereby $*$-conjugation is implemented
in braided geometry by reversal of braid crossings. The equivalence of the
action and conjugate action via the
braided-parity operator also applies to all orders of products of $q$-spacetime
generators. This is the sense within braided geometry  in which the operators
$c_i$ are `antihermitian'.  This also holds for the momentum $p^i$ generators
as the main result in \cite{Ma:qsta}.

The present work suggests the possibility of a systematic theory of massless
spinning particles based on invariance under the $q$-conformal group.
This will be attempted elsewhere. Classically, it requires the construction of
fields with conformal weights defined as sections of certain vector bundles
over compactified spacetime. In the $q$-deformed case one needs therefore
nontrivial quantum homogeneous spaces and their associated bundles, for example
along the lines in \cite{BrzMa:gau}.

\appendix
\section{Abstract Results}

Most of the formulae for the $q$-conformal group in the text above have been
given at the level of R-matrices and matrix relations. In principle, one also
has to check a large number of non-quadratic relations, in particular
associated with the $\vecl^\pm$ generators (they are not independent). These
are needed to form a (quasitriangular) Hopf algebra in the Lorentz sector.
Fortunately, such details are ensured by the abstract braided group and quantum
group constructions underlying the R-matrix formulae. This is given for the the
extended $q$-Poincar\'e algebra in \cite{Ma:qsta} and we extend this now for
the $q$-conformal case. The basis for the latter is a recent construction
\cite{Ma:con} which associates to a braided group $B$ in the category of
$H$-modules ($H$ a quasitriangular Hopf algebra), a new quasitriangular Hopf
algebra built from $B,H,B^*$, called the {\em double-bosonisation} of $B$. Here
we state without proof the relevant left-module version of the
double-bosonisation formulae (different  right-module conventions are used in
\cite{Ma:con}, for the purposes there). Then we study $*$-structures in this
abstract setting, which is the new result of this appendix.

Familiarity with abstract quantum group\cite{Ma:book} and braided
group\cite{Ma:introp} techniques is assumed. In particular, $\Delta h=h\o\tens
h\t$ denotes the coproduct of $h\in H$ and $\und\Delta b=b\Bo\tens b\Bt$ the
braided coproduct of $b\in B$. Also, $\CR=\CR\uo\tens\CR\ut$ denotes the
quasitriangular structure of $H$, and $\cv=\CR\uo S\CR\ut$.

Let $B$ be a braided group with invertible braided antipode in the braided
category of left $H$-modules, dual to another braided group $C$. So
$B=C^\star$. (In the infinite-dimensional case, we suppose a duality pairing
$\ev:B\tens C\to k$ of braided groups.)  The {\em double bosonisation} of $B$
is the Hopf algebra $U(B)$ containing $B, C^{\rm op}, H$ as subalgebras and the
cross relations, coproduct and antipode cf\cite{Ma:con}
\ceqn{dbos}{ hb=(h\o\la b)h\t,\quad hc=(h\t\la c)h\o\\
b\Bo\CR\ut c\Bo\ev(\CR\uo\la b\Bt,c\Bt)=\ev(b\Bo,\CR\ut\la c\Bo) c\Bt \CR\uo
b\Bt\\
\Delta b=b\Bo \CR\ut\tens \CR\uo\la b\Bt,\quad \Delta c=\CR\ut\la c\Bo\tens
c\Bt\CR\uo\\
Sb=(\cu\CR\uo\la \und Sb) S\CR\ut,\quad Sc=\CR\umo \und
S^{-1}(\cv^{-1}\CR\umt\la c)}
where $\la$ denotes the action of $H$ whereby $B,C$ live in the braided
category of $H$-modules. The pairing is assumed covariant, so $\ev(h\la
b,c)=\ev(b,(Sh)\la c)$. The unit and counit are the trivial tensor product ones
and $H$ has its usual coproduct and antipode (it is a sub-Hopf algebra).
Similar proofs to those in \cite{Ma:con} show that this defines a Hopf algebra.
The bosonisation $B\lbiprod H$ appears as a sub-Hopf algebra and a certain
`conjugate bosonisation' generated by $C^{\rm op}, H$ also appears as a
sub-Hopf algebra.

When the pairing is non-degenerate, we have as a formal power series a
canonical
element $\exp=e_\alpha\tens f^\alpha$ for the pairing, where $\{e_\alpha\}$ is
a basis of $C$ and $\{f^\alpha\}$ a dual basis. Its inverse in the algebra
$C^{\rm op}\tens B$ is $\exp^{-1}=(\und S e_\alpha\tens f^\alpha)$ from the
pairing axioms. In this case the double-bosonisation is quasitriangular with
\eqn{dbosR}{\CR_{U(B)}=\CR\exp^{-1},}
where we view $\CR$ and $\exp$ in $U(B)\tens U(B)$.

This describes left-handed version of the formulae in \cite{Ma:con}. It
underlies the formulae in Sections~2: we take $H$ generated by
$\vecl^\pm,\dila$ and $B$ generated by $\vecp$. We take $C$ generated by
$\vecc$ which are dual to the $\vecp$ in the usual way except scaled so that
$\ev(p^i,c_j)=(q-q^{-1})^{-1}\delta^i{}_j$. We then use the same methods as in
\cite{Ma:poi}\cite{Ma:qsta} for the calculation of the extended $q$-Poincar\'e
algebra as $B\lbiprod H$. Similarly for the conjugate bosonisation generated by
$H,C^{\rm op}$. The remaining cross relations are
$[p^i,c_j]=\ev(p^i,c_a)\<t^a{}_j g,\CR\ut\>\CR\uo-\CR\ut\<St^i{}_a
g^{-1},\CR\uo\>\ev(p^a,c_j)$ from (\ref{dbos}) and the linear form of the
braided coproducts on $B,C$. It is convenient to compute the action of $\CR$
here on $C$ as evaluation against the coaction of the matrix quantum group dual
to $H$, with generators $\vect,g$. We derive the formula in Proposition~2.1 in
this way. The same method gives the spinorial formulae in Section~4.

Also given in \cite{Ma:con} is a fundamental representation of $U(B)$. In our
conventions it appears on $C$, making it into a left $U(B)$-module algebra as
follows: first, $B$ acts on $C$ by the braided left coregular representation
studied in \cite[Sec. 2]{Ma:qsta}. Together with the given action $\la$ of $H$
on $C$, we have\cite[Cor.~2.2]{Ma:qsta} a covariant action of $B\lbiprod H$ on
$C$ (it defines  the action of the extended $q$-Poincar\'e algebra on
spacetime). Secondly, $C$ acts on itself by the right braided-adjoint action.
This is given by the diagrams in \cite{Ma:lie} reflected in a mirror followed
by reversal of braid crossings. We view this as a left action of $C^{\rm op}$
on itself. Cf\cite{Ma:bos}, these actions fit together to give an action of
$U(B)$ covariantly on $C$ (i.e. respecting its product). Explicitly,
\eqn{dbosact}{ b\la x=\ev(\und S^{-1}b, x\opBo)x\opBt,\quad c\la x=(\CR\ut\la
\und S c\Bo)(\CR\uo\la x)c\Bt}
when acting on $x\in C$. Here $x\opBo\tens x\opBt=\CR\umo\la x\Bt\tens
\CR\umt\la x\Bo$ is the opposite braided coproduct of $C$. This action
underlies the formulae in Section~5. The action of $p^i\in B$ is by
braided-differentiation as studied in \cite{Ma:qsta}. The action of $c_i\in C$
is computed as
\[c_i\la x_j=(\CR\ut\la \und S c_i) x_a \<\CR\uo,t^a{}_j g\>+x_jc_i
= x_j c_i -(\dila^{-1} Sl^-{}^a{}_j\la c_i)x_a={x_jx_i-x_ax_b
R^b{}_j{}^a{}_i\over q-q^{-1}}\]
where $x_i\in C$ are the usual (not scaled) generators dual to $p^i$ (so
$c_i=(q-q^{-1})x_i$), and $\und S x_i=-x_i$. This derives the action used in
Section~5.

Next we move on to new abstract considerations beyond  \cite{Ma:con}. We
suppose that $B,C$ are $*$-braided groups in the usual sense\cite{Ma:star}, $H$
is a real-quasitriangular Hopf $*$-algebra and its action on $B$ is unitary in
the Hopf algebraic sense. Thus
\eqn{bosuni}{ (h\la b)^*=(Sh)^*\la b^*,\quad (h\la c)^\star=(S(h^*))\la
c^\star.}
As explained in \cite{Ma:qsta}, the second formula is dictated by the first one
 and braided group duality. We use $\star$ to define the $*$-structure on
$C^{\rm op}$ as the same antilinear map.

\begin{propos} In this setting, the double bosonisation $U(B)$ is a quasi-$*$
Hopf algebra with cocycle $\CR$ viewed in $U(B)\tens U(B)$. Moreover,
\[ (*\tens *)\circ\Delta\circ *=\exp^{-1}(\Delta\ )\exp\]
in $U(B)$.
\end{propos}
\proof It is proven in \cite{Ma:qsta} that, in this setting, $B\lbiprod H$
becomes a quasi-$*$ Hopf algebra with cocycle $\CR$. By a similar calculation,
we find that the conjugate bosonisation generated by $C^{\rm op}, H$ is also a
quasi-* Hopf algebra with the same cocycle $\CR$. We verify that these
$*$-structures are compatible with the cross-relations in (\ref{dbos}).
Applying $*$ to both sides:
\align{&&\equad\equad (b\Bo\CR\ut c\Bo)^*\overline{\ev(\CR\uo\la
b\Bt,c\Bt)}=c\Bo{}^*\CR\ut{}^* b\Bo{}^*\ev((\CR\uo\la b\Bt)^*,c\Bt{}^*)\\
&=&c^*\Bt \CR\uo b^*\Bt\ev(S^{-1}\CR\ut\la b^*\Bo,c^*\Bo)=c^*\Bt \CR\uo
b^*\Bt\ev(b^*\Bo,\CR\ut\la c^*\Bo)\\
&=&\ev(\CR\uo\la b^*\Bt,c^*\Bt) b^*\Bo\CR\ut c^*\Bo=\ev( b^*\Bt,S\CR\uo\la
c^*\Bt)b^*\Bo\CR\ut c^*\Bo\\
&=&\ev(b\Bo{}^*,(\CR\ut\la c\Bo)^*)b\Bt{}^* \CR\uo{}^* c\Bt{}^*
=\overline{\ev(b\Bo,\CR\ut\la c\Bo)} (c\Bt \CR\uo b\Bt)^*}
using $\overline{\ev(b\tens c)}=\ev(b^*,c^*)$, reality of $\CR$ in the sense
$\CR^{*\tens *}=\CR_{21}$, our assumption (\ref{bosuni}), invariance of $\ev$
and the cross-relations in (\ref{dbos}) applied to $b^*,c^*$. This checks
consistency of the relations under $*$ and implies that we have a $*$-algebra
structure on $U(B)$.  Since its two sub-Hopf algebras mentioned above are
quasi-$*$ Hopf algebras with cocycle $\CR$, it becomes a quasi-$*$ Hopf algebra
 as well, with the same cocycle.

Since $U(B)$ is also (in the non-degenerately paired case) quasitriangular via
(\ref{dbosR}), we deduce that $(*\tens *)\circ\Delta\circ
*=\CR^{-1}(\tau\circ\Delta\ )\CR=\CR^{-1}\CR_{U(B)}(\Delta\
)\CR_{U(B)}^{-1}\CR=\exp^{-1}(\Delta\ )\exp$, as stated. \endproof

We see from this proposition that the plane wave $\exp$ controls the extent
that the double bosonisation fails to be a Hopf $*$-algebra in the usual sense.
This, in turn, expresses the sense in which the tensor product of unitaries
fails to be unitary: they are unitary only up to a cocycle isomorphism
expressed
by the action of $\exp$.

{}From the theory of quasi-$*$ Hopf algebras in \cite[Lemma~4.7]{Ma:qsta} it is
known that if a quasi-$*$ Hopf algebra $\CH$ acts covariantly on a $*$-algebra
$C$ by $\la$, then the conjugate quasi-$*$ Hopf algebra (with coproduct
$\bar\Delta=(*\tens *)\circ\Delta\circ *$) acts covariantly on $C$ by a
conjugate action $\bar\la$ defined by
\eqn{conjaction}{h\bar\la x=(S(h^*)\la x^*)^*}
for all $h\in \CH$ and $x\in C$.

\begin{theorem} The conjugate of the action of $U(B)$ as a quasi-$*$ Hopf
algebra acting on $C$ is $h\bar\la x=h\la x$ and $b\bar\la x=\ev(b,x\Bo)x\Bt$
as in \cite{Ma:qsta}, and
\[ c\bar\la x=(\CR\umo\la \und S^{-1} c\opBo)(\CR\umt\la x)c\opBt.\]
Moreover, $(\ )\la\und S x=\und S((\ )\bar\la x)$, i.e. the action and
conjugate action of $U(B)$ are intertwined by the braided antipode of $C$.
\end{theorem}
\proof The conjugate actions of $h\in H$ and $b\in B$ are covered in the
conjunction of \cite[Cor.~2.4]{Ma:qsta} and \cite[Prop.~4.8]{Ma:qsta}. To this
we add now the conjugate of the action of $c\in C$. We compute:
\align{c\bar\la x\equad&&=(S(c^*)\la
x^*)^*=\left(\CR\umo\la((\cv^{-1}\CR\umt\la\und S^{-1}c^*)\la
x^*)\right)^*=S(\CR\umo{}^*)\la \left( (\cv^{-1}\CR\umt\la \und S^{-1} c^*)\la
x^*\right)^*\\
&&=S(\CR\umo{}^*)\la\left((\CR\ut\la\und S(\cv^{-1}\CR\umt\la\und S^{-1}
c^*)\Bo)(\CR\uo\la x^*)(\cv^{-1}\CR\umt\la\und S^{-1} c^*)\Bt\right)^*\\
&&=S(\CR\umo{}^*)\la\left((\cv^{-1}\CR\umt\la\und S^{-1} c^*)^*\Bo (
S(\CR\uo{}^*)\la x)( S(\CR\ut{}^*)\la\und S (\cv^{-1}\CR\umt\la\und S^{-1}
c^*)^*\Bt)\right)\\
&&=\CR\umt\la\left((\cu^{-1}\CR\umo\la\und S^{-1}c)\Bo (\CR\ut\la
x)(\CR\uo\la\und S(\cu^{-1}\CR\umo\la\und S^{-1}c)\Bt)\right)\\
&&=\CR_2\umt\CR_1\umt\la\left((\cu^{-1}\CR_1\ut\CR_2\uo\CR_1\umo\CR\umo\la\und
S^{-1} c\Bt)(\CR\ut\la x)(\CR\uo\cu^{-1}\CR_1\uo\CR_2\ut\CR_2\umo\CR\umt\la
c\Bo)\right)}
where we use the antipode of $U(B)$ from (\ref{dbos}) and repeatedly use
(\ref{bosuni}). For the fifth equality we use the axiom $c\Bo{}^*\tens
c\Bt{}^*=c^*\Bt\tens c^*\Bo$ for $*$-braided groups.  We then use the reality
property of $\CR$, which also implies that $\cv^{-1}{}^*=\cv^{-1}$. Here
$S\cv^{-1}=\cu^{-1}=\CR\ut S^2\CR\uo$. For the last equality we use covariance
of $\und \Delta$ under the action of $H$ (along with standard facts about
quasitriangular Hopf algebras to compute $\Delta\cu^{-1}$ and
$(\Delta\tens\id)\CR^{-1}$),  and the braided anticomultiplicativity
$\und\Delta\und S^{-1}c=\CR\umo\la\und S^{-1}c\Bt\tens \CR\umt\la\und
S^{-1}c\Bo$   from \cite{Ma:introp}. Numerical suffices on $\CR,\CR^{-1}$ are
used to distinguish the various copies. The remaining steps are a tedious but
straightforward computation: we use the QYBE for $\CR$ to cancel some of the
$\CR$ factors. Then we compute the action of $\CR_1\umt\CR_2\umt$   on products
 using covariance, converting coproducts on $\CR^{-1}$ into more copies of
$\CR^{-1}$. Using $\cu(\ )\cu^{-1}=S^2$ and $\CR\umt\cu^{-1}\CR\umo=1$, we can
then cancel most of the $\CR,\CR^{-1}$ factors to obtain the result stated.
Here   $c\opBo\tens c\opBt$ denotes the braided opposite coproduct of $C$ as
usual.

Next, we show that the action and conjugate action are intertwined by the
braided antipode $\und S$ of the copy of $C$ in which we are acting. For the
action of $h\in H$ this is covariance of the braided antipode. For $b\in B$
this is \cite[Cor.~2.4]{Ma:qsta}. To this we now add:
\align{\und S(c\bar \la x)\equad&&= (\CR\ut\la \und S
c\opBt)\left(\CR\uo\la\und S( (\CR\umo\la\und S^{-1}c\opBo)(\CR\umt\la
x))\right)\\
&&=(\CR_1\ut\CR_2\ut\la\und S c\opBt)(\CR_1\uo\la\und S x)(\CR_2\uo\la
c\opBo)=(\CR\ut\la\und S c\Bo)(\CR\uo\la \und S x)c\Bt=c\la\und S x}
using braided-antimultiplicativity of $\und S$ twice. This part of the proof
can also be done diagrammatically.  \endproof

This therefore extends the abstract unitarity and quasi-$*$ considerations for
bosonisations and the extended $q$-Poincar\'e algebra in \cite{Ma:qsta} to
double-bosonisations and the $q$-conformal algebra. It is used in Sections~3
and~5.

%\bibliographystyle{unsrt}
%\bibliography{biblio}

\end{document}